\documentclass{amsart}

\newcommand{\Exp}[2]{\ensuremath{\exp_{#1}\!{#2}}}
\newcommand{\U}{\mathcal{U}}
\newcommand{\naturalnumbers}{\ensuremath{\mathbf{N}}}
\def\co{\colon\thinspace}

\newtheorem{theorem}{Theorem}
\newtheorem{lemma}{Lemma}
\newtheorem{corollary}{Corollary}
\newtheorem*{conjecture}{Conjecture}

\begin{document}

\title{Connectivity of finite subset spaces of cell complexes}
\author{Christopher Tuffley}
\address{Department of Mathematics, University of California at Davis \\
         One Shields Avenue, Davis, CA 95616-8633, U.S.A.}
\email{tuffley@math.ucdavis.edu}

\subjclass{55R80 (54B20 55Q52)}
\keywords{Configuration spaces, finite subset spaces, cell complexes,
homotopy groups}

\begin{abstract}
The $k$th finite subset space of a topological space $X$ is the space
\Exp{k}{X}\ of non-empty subsets of $X$ of size at most $k$,
topologised as a quotient of $X^k$. Using results from our earlier
paper on the finite subset spaces of connected graphs we show that 
the $k$th finite subset space of a connected cell complex is
$(k-2)$--connected, and $(k-1)$--connected if in addition the underlying
space is simply connected. We expect \Exp{k}{X}\ to be $(k+m-2)$--connected
if $X$ is an $m$--connected cell complex, and reduce proving 
this to the problem of proving it for finite wedges of $(m+1)$--spheres.
Our results complement a
theorem due to Handel that for path-connected Hausdorff $X$ the map on
$\pi_i$ induced by the inclusion
$\Exp{k}{X}\hookrightarrow\Exp{2k+1}{X}$ is zero for all $k\geq 1$ and
$i\geq 0$. 
\end{abstract}

\maketitle

\subsection{Introduction}

The $k$th finite subset space of a topological space 
$X$ is the space $\Exp{k}{X}$ of non-empty subsets of $X$
of size at most $k$, topologised as a quotient of $X^k$ via the map
\[
(x_1,\ldots,x_k)\mapsto\{x_1\}\cup\cdots\cup\{x_k\}.
\]
The construction is functorial, with $f\co X\rightarrow Y$
inducing $\Exp{k}{f}\co\Exp{k}{X}\rightarrow\Exp{k}{Y}$ sending
$\alpha\subseteq X$ to $f(\alpha)\subseteq Y$, and moreover \Exp{k}{f}
and \Exp{k}{g} are homotopic if $f$ and $g$ are, so that $\Exp{k}{}$
is in fact a homotopy functor.

Handel~\cite{handel00} has shown that for path-connected Hausdorff $X$
the map on $\pi_i$ induced by the inclusion 
$\Exp{k}{X}\hookrightarrow\Exp{2k+1}{X}$ is zero. Using results from
our paper~\cite{graphs03} on the finite subset spaces of connected
graphs we complement this
result, proving the following connectivity theorem for the 
finite subset spaces of a connected cell
complex:

\begin{theorem}
\label{vanish.th}
The $k$th finite subset space of a connected cell complex $X$ is
$(k-2)$--connected, in other words 
$\pi_i(\Exp{k}{X})$ vanishes for $i\leq k-2$.
\end{theorem}

Using results from~\cite{surfaces03} the conclusion may be strengthened
to \Exp{k}{X}\ is $(k-1)$--connected if $X$ is a simply connected 
cell complex, and we expect further that 
\Exp{k}{X}\ should be $(k+m-2)$--connected if $X$ is an 
$m$--connected cell complex. We discuss this in 
section~\ref{vanishdiscuss.sec} after proving the result in 
section~\ref{vanishproof.sec}.

\subsection{Proof of Theorem~\ref{vanish.th}}
\label{vanishproof.sec}

We note that Theorem~\ref{vanish.th} is immediate for $k=2$, since 
the finite subset spaces
of a path-connected space are again path-connected, so in proving it
we assume $k\geq 3$. In this case the conclusion is equivalent to the
statement that \Exp{k}{X}\ is simply connected with vanishing
reduced homology in dimensions less than $k-1$, by the Hurewicz
theorem, and we shall freely use this formulation. 
We first prove the 
result for finite complexes and then pass to the infinite case via
a compactness argument.

\begin{proof}[Proof of Theorem~\ref{vanish.th} for finite $X$]
The proof is by induction on the dimension of $X$, with the base
case $\dim X=1$ given by our paper~\cite{graphs03}, in which we studied
the finite subset spaces of connected finite graphs.
The keys to the inductive step are the following two observations:
\begin{enumerate}
\item
\label{cover.obs}
If $P_1,\ldots,P_{k+1}$ are disjoint subsets of $X$ then a $k$ element
subset of $X$ must lie in $X\setminus P_\ell$ for some $\ell$.
\item
\label{n-1.obs}
If $X$ is a connected finite $(n+1)$--complex, $n\geq 1$, and $P\subseteq X$ 
lies in the open $(n+1)$--cells and intersects each in a non-empty
finite set, then $X\setminus P$ has the homotopy type of a connected finite
$n$--complex.
\end{enumerate}
We use these in conjunction with the following lemma, itself
proved inductively using the Mayer-Vietoris sequence.
\begin{lemma}
\label{mv.lem}
Let $Y$ be a union of open sets $U_1,\ldots,U_r$ such that 
$U_1\cap\cdots\cap U_r$ is non-empty and each 
$U_{i_1}\cap\cdots\cap U_{i_s}$ has vanishing reduced homology in 
dimensions less than $j$. Then $Y$ has vanishing reduced
homology in dimensions less than $j$ also.
\end{lemma}

Suppose the theorem holds for connected finite $n$--complexes, some
$n\geq 1$, and let $X$ be a connected finite $(n+1)$--complex.
Let $v$ be a vertex of $X$ and let $P_1,\ldots,P_{k+1}$ be disjoint
subsets of $X$ each consisting of exactly one point from each open 
$(n+1)$--cell. By observation~(\ref{cover.obs}) the sets
\[
\U_\ell = \Exp{k}{(X\setminus P_\ell)}
\]
cover \Exp{k}{X}, and each $\U_\ell$ is open since $(X\setminus P_\ell)^k$
is open in $X^k$. Moreover the intersection $\U_1\cap\cdots\cap\U_{k+1}$ 
contains $\{v\}$ and is therefore non-empty.

Consider 
\[
\U_{\ell_1}\cap\cdots\cap \U_{\ell_s} = 
  \Exp{k}{(X\setminus (P_{\ell_1}\cup\cdots\cup P_{\ell_s}))}.
\]
By observation~(\ref{n-1.obs}) the space
$X\setminus (P_{\ell_1}\cup\cdots\cup P_{\ell_s})$ has the homotopy
type of a connected finite $n$--complex, so by the inductive hypothesis
$\U_{\ell_1}\cap\cdots\cap \U_{\ell_s}$ has vanishing reduced homology
in dimensions less than $k-1$. It follows that the hypotheses of
Lemma~\ref{mv.lem} are satisfied by the cover $\{\U_\ell\}$ with
$j=k-1$, and we conclude that $\Exp{k}{X}$ has vanishing reduced
homology in dimensions less than $k-1$.

To complete the inductive step it remains to show that $\Exp{k}{X}$ is
simply connected. This follows immediately from the van Kampen theorem
applied to the cover $\{\U_\ell\}$ with basepoint $\{v\}$. By the
inductive hypothesis $\pi_1(\U_\ell,\{v\})\cong\{1\}$ for all $\ell$, and
since each $\U_j\cap\U_\ell$ is path-connected we have immediately
that $\pi_1(\Exp{k}{X},\{v\})\cong\{1\}$ also. 
\end{proof}

To pass to the infinite dimensional case we use the following lemma, which
we prove with no assumptions on $X$.

\begin{lemma}
If $C\subseteq \Exp{k}{X}$ is
compact in \Exp{k}{X}\ then 
$\bigcup C = \bigcup_{\alpha\in C} \alpha$
is compact in $X$.
\end{lemma}

\begin{proof}
Given an open cover $\mathcal O$ of $\bigcup C$ the set
\[
\mathcal{O}' = 
   \bigl\{\Exp{k}{\textstyle\bigcup O} \big| 
                          O\subseteq\mathcal{O}\mbox{ finite}\bigr\}
\]
is an open cover of $C$. Extracting a finite subcover 
$\bigl\{\Exp{k}{\bigcup O_1},\ldots,\Exp{k}{\bigcup O_m}\bigr\}$ from 
$\mathcal{O}'$ we obtain a finite
subcover $O_1\cup\cdots\cup O_m$ of $\mathcal{O}$.
\end{proof}

Since compact subspaces of cell complexes lie in finite subcomplexes
we have immediately:

\begin{corollary}
If $X$ is a cell complex and $C\subseteq \Exp{k}{X}$ is compact, then
$C\subseteq\Exp{k}{A}$ for some finite subcomplex $A$ of $X$.
\label{insubcomplex.cor}
\end{corollary}

\begin{proof}[Proof of Theorem~\ref{vanish.th} for $X$ infinite]
Let $[\phi]\in\pi_i(\Exp{k}{X},\{v\})$ for some vertex $v$ of $X$ and
$i\leq k-2$. By Corollary~\ref{insubcomplex.cor}
$\phi(S^i)$ lies in \Exp{k}{A}\ for
some finite subcomplex $A$ of $X$, and we may take $A$ to be connected
since $X$ is.
By the finite case of the theorem $[\phi]$ is trivial in 
$\pi_i(\Exp{k}{A},\{v\})$, and so in $\pi_i(\Exp{k}{X},\{v\})$ also.
\end{proof}

\subsection{Discussion}
\label{vanishdiscuss.sec}

Theorem~\ref{vanish.th} is consistent
with the following conjecture on the finite subset spaces of cell complexes, 
as the theorem follows from the conjecture together with Handel's inclusion
result. 
We restrict our attention to
complexes with a single vertex in each component, with no 
loss of generality up to homotopy.

\begin{conjecture}
Let $X$ be a finite $n$--complex with $c$ components, each 
containing a single vertex. Then $\Exp{k}{X}$ has a
cell structure obtained from $\Exp{k-1}{X}$
by adding cells of dimensions $k-c\leq i\leq nk$.
\end{conjecture}

The conjecture is true in the connected case for $n=1$
by Lemma~1 of~\cite{graphs03}, and for $n=2$ by 
Theorem~6 of~\cite{surfaces03}. 
To see that it implies Theorem~\ref{vanish.th} note that it implies
that the homotopy groups of $\Exp{k}{X}$ stabilise as $k$ increases.
By Handel's result the stable groups must be zero when $X$ is
connected, and careful attention to the point at which the
stabilisation occurs gives the bound in the theorem.

Jacob Mostovoy (private communication) has indicated that finite
subset spaces of cell complexes may be shown to have cell structures
using the machinery of simplicial sets, which is described
in May's book~\cite{may67} or Curtis's article~\cite{curtis71}. 
Given a simplicial set $K$ we let \Exp{j}{K}\ be the simplicial
set whose $n$--simplices are subsets of size at most $j$ of the
$n$--simplices of $K$, and whose face and degeneracy operators
are the face and degeneracy operators of $K$ acting elementwise. 
Then if $X$ is the geometric realisation of $K$, \Exp{j}{X}\ will
be the geometric realisation of \Exp{j}{K}, showing that
triangulated spaces have triangulated finite subset spaces. However,
the triangulations obtained in this way 
apparently do not satisfy the lower bound $k-c$ needed to prove
Theorem~\ref{vanish.th}.  The bound comes from the following
heuristic, which is motivated by the form of the lexicographic cell
structures of~\cite{surfaces03}. 
We suppose that \Exp{k}{X}\ has a cell structure such that
for each open cell $e$ of $X$ the set map
\[
\Exp{k}{X}\rightarrow\naturalnumbers : \Lambda \mapsto |\Lambda\cap e|
\]
is constant on each open cell $\sigma$ of \Exp{k}{X}, and moreover
that the dimension of $\sigma$ is at least the common number
of points in $X$ less the $0$--skeleton for $\Lambda\in\sigma$.
In particular we suppose that the vertices of \Exp{k}{X}\ can be
chosen to be subsets of the vertices of $X$, without adding more
through subdivision.
If $X$ has a single vertex in each component then an open cell of 
$\Exp{k}{X}\setminus\Exp{k-1}{X}$ in such a cell structure 
would have dimension at least $k-c$.

Theorem~\ref{vanish.th} can be strengthened for simply 
connected complexes, and we expect that it can be strengthened further
for $m$--connected cell complexes.  With this in mind we prove
the following theorem, showing it suffices to prove any strengthened 
result for wedges of spheres. 

\begin{theorem}
Suppose that finite wedges of $(m+1)$--spheres have $r$--connected
$k$th finite subset spaces. Then $m$--connected cell complexes have
$r$--connected $k$th finite subset spaces also.
\label{wedgesofspheres.th} 
\end{theorem}

\begin{proof}
We simply adapt the inductive step of the proof of 
Theorem~\ref{vanish.th}, replacing
observation~(\ref{n-1.obs}) with the following observation~(\ref{n-1.obs}$')$:
\begin{enumerate}
\item[\ref{n-1.obs}$'$.]
If $X$ is an $m$--connected finite $(n+1)$--complex, $n\geq m+1$, 
and $P\subseteq X$ 
lies in the open $(n+1)$--cells and intersects each in a non-empty
finite set, then $X\setminus P$ has the homotopy type of an $m$--connected 
finite $n$--complex.
\end{enumerate}
If $X$ is an $m$--connected finite cell complex then up to homotopy
we may assume that the $(m+1)$--skeleton of $X$ is a finite wedge
of $(m+1)$--spheres. The base for the induction is then given by 
hypothesis and the argument goes through exactly as before.
\end{proof}

By~\cite[Theorem 6]{surfaces03} $\Exp{k}{\bigvee_n S^2}$ has a 
cell structure obtained from $\Exp{k-1}{\bigvee_n S^2}$ by adding
cells in dimensions $k$ and higher, and using Handel's inclusion
result and Theorem~\ref{wedgesofspheres.th} we conclude that 
simply connected cell complexes have $(k-1)$--connected $k$th finite
subset spaces. More generally, the construction outlined 
in~\cite[section 2.4]{surfaces03} should yield cell structures for
the finite subset spaces of wedges of spheres 
in which $\Exp{k}{\bigvee_n S^{m+1}}$
is obtained from $\Exp{k-1}{\bigvee_n S^{m+1}}$ by adding cells in 
dimensions $k+m-1$ and higher. Verifying the details of this construction
would show that $m$--connected cell complexes have $(k+m-2)$--connected
$k$th finite subset spaces.

We conclude with an example showing  
the necessity of the connectedness hypothesis in Theorem~\ref{vanish.th}.
Consider the third finite subset space of a
pair of circles, $\Exp{3}{(S^1\amalg S^1)}$. This has three
connected components, two ``pure'' components consisting of 
subsets contained entirely in one or the other component circle,
and a ``mixed'' component consisting of subsets meeting both. 
The two pure components are each homeomorphic to $S^3$
(Bott~\cite{bott52}; see also~\cite{circles02}),
but the mixed component is formed by gluing two copies of 
$\Exp{2}{S^1}\times \Exp{1}{S^1}\cong \mbox{M\"ob}\times S^1$ along
their boundary. The gluing interchanges the roles of the boundary
of the M\"obius strip 
and the $S^1$ direction, and $\pi_1$ of the resulting three manifold
has presentation $\langle a,b| [a,b^2]=[a^2,b]=1\rangle$.

\end{document}